\documentclass[11pt, twoside, eqno]{article}
\usepackage{latexsym}
\usepackage{amssymb}
\usepackage{amsfonts}
\begin{document}
\begin{center}

\noindent {\bf \Large Remark on quasi-ideals of ordered
semigroups}\bigskip

\medskip

{\bf Niovi Kehayopulu}\bigskip

{\small Department of Mathematics,
University of Athens \\
15784 Panepistimiopolis, Athens, Greece \\
nkehayop@math.uoa.gr }
\end{center}
\date{ }

\begin{center}
{\bf Abstract} \end{center} {\small The aim is to correct part of the  
Remark 3 of my paper ``On regular, intra-regular ordered semigroups" 
in Pure Math. Appl. (PU.M.A.) 4, no. 4 (1993), 447--461. On this 
occasion, some further results and the similarity between the 
$po$-semigroups and the $le$-semigroups is discussed.\bigskip

\medskip

}
\bigskip

\section {Introduction and prerequisites}
There is a typing mistake in Remark 3 in [8], corrected by hand in 
the reprints: On page 458, line --4, the set $X\cup {\Big( }(XS]\cap 
(SX]{\Big )}$
should be replaced by ${\bigg (}X\cup {\Big( }(XS]\cap (SX]{\Big 
)}{\bigg ]}$. In this Remark we have seen that, for the subset 
$Q:=X\cup {\Big( }(XS]\cap (SX]{\Big )}$ of $S$, we have

1) $(QS]\cap (SQ]\subseteq Q$ and

2) if $T$ is a quasi-ideal of $S$ such that $X\subseteq T$, then 
$Q\subseteq T$, thus $(Q]\subseteq (T]=T$. In addition, we have 
$((Q]]=Q$. As a consequence, the set $(Q]$ is the quasi-ideal of $S$ 
generated by $X$. Indeed: $${\Big(}(Q]S{\Big]}\cap 
{\Big(}S(Q]{\Big]}={\Big(}(Q](S]{\Big]}\cap 
{\Big(}(S](Q]{\Big]}\subseteq (QS]\cap (SQ]\subseteq Q.$$

Though the following does not affect the proof, on page 459, line 13, 
it is better to write ${\bigg (}X\cup {\Big( }(XS]\cap (SX]{\Big 
)}{\bigg ]}$ instead of ${\bigg (}X\cup {\Big( }(XS]\cap (SX]{\Big 
)}{\bigg )}$; furthermore in line 16 of the same page $X^2S$ should 
be replaced by $(X^2S]$ taking it into account in the rest of the 
proof of the Remark. For convenience, we give here detailed proofs of 
our arguments. \medskip

For an ordered semigroup $S$ and a subset $H$ of $S$, we denote by 
$(H]$ the subset of $S$ defined by $$(H]:=\{t\in S \mid t\le h \mbox 
{ for some } h\in H\}.$$A nonempty subset $A$ of $S$ is called a 
quasi-ideal of $S$ if (1) $(AS]\cap (SA]\subseteq A$ and (2) if $a\in 
A$ and $S\ni b\le a$, then $b\in A$ (equivalently $(A]\subseteq A$, 
which in turn is equivalent to $(A]=A$). It is called a bi-ideal of 
$S$ if (1) $ASA\subseteq A$ and (2) if $a\in A$ and $S\ni b\le a$, 
then $b\in A$. Every quasi-ideal of $S$ is a bi-ideal of $S$ as well. 
Indeed, if $B$ is a quasi-ideal of $S$, then $BSB\subseteq (BS]\cap 
(SB]\subseteq B$. A nonempty subset $A$ of $S$ is called a left 
(resp. right) ideal of $S$ if (1) $SA\subseteq A$ (resp. $AS\subseteq 
A$) and (2) if $a\in A$ and $S\ni b\le a$, then $b\in A$. For a 
subset $A$ of $S$, the set $(A\cup AS]$ is the right ideal of $S$ 
generated by $A$, and the set $(A\cup SA]$ is the left ideal of $S$ 
generated by $A$. An ordered semigroup $S$ is called intra-regular if 
for every $a\in S$ there exist $x,y\in S$ such that $a\le xa^2y$. 
This is  equivalent to saying that $a\in (Sa^2S]$ for every $a\in S$ 
or $A\subseteq (SA^2S]$ for every $A\subseteq S$ (cf., for example 
[7]).

We mention the properties we use in the paper: Clearly $S=(S]$, and 
for subsets $A, B$ of $S$, we have the following:
if $A\subseteq B$, then $(A]\subseteq (B]$; $A\subseteq (A]$; 
${\Big(}(A]{\Big]}=(A]$; $(A](B]\subseteq (AB]$; 
${\Big(}(A](B]{\Big]}=(AB]$ (cf., for example [6]).

\section{Main results}

If $S$ is an ordered semigroup, for a subset $X$ of $S$, the set 
$$Q:={\Bigg (}X\cup {\Big( }(XS]\cap (SX]{\Big )}{\Bigg ]}$$is the 
quasi-ideal of $S$ generated by $X$. It is mentioned without proof in 
[9]. In spite of whatever we already said at the beginning of the 
introduction which gives a complete proof of our argument, we think 
it is interesting to give an independent detailed proof which is the 
following:

1) $(QS]\cap (SQ]\subseteq Q$. Indeed:
\begin{eqnarray*}QS&=&{\Bigg (}X\cup {\Big( }(XS]\cap (SX]{\Big 
)}{\Bigg ]}S={\Bigg (}X\cup {\Big( }(XS]\cap (SX]{\Big )}{\Bigg 
]}(S]\\&\subseteq& {\Big (}X\cup (XS]{\Big ]}(S]\subseteq 
{\Bigg(}{\Big(}X\cup (XS]{\Bigg)}S{\Bigg]}\\&=&{\Big(}XS\cup 
(XS]S{\Big]}.
\end{eqnarray*}Since $(XS]S=(XS](S]\subseteq (XS^2]\subseteq (XS]$, 
we have$$QS\subseteq {\Big(}XS\cup 
(XS]{\Big]}={\Big(}(XS]{\Big]}=(XS],$$so $(QS]\subseteq 
{\Big(}(XS]{\Big]}=(XS]$. Similarly we get so $(SQ]\subseteq (SX]$, 
thus we have\begin{eqnarray*}(QS]\cap (SQ]&\subseteq&(XS]\cap 
(SX]\subseteq X\cup {\Big(}(XS]\cap (SX]{\Big)}\\&\subseteq & 
{\Bigg(}X\cup {\Big(}(XS]\cap (SX]{\Big)}{\Bigg ]}=Q.\end{eqnarray*}

2) If $a\in Q$ and $S\ni b\le a$, then $b\in Q$. In fact:\\
Since $a\in Q:={\Bigg (}X\cup {\Big( }(XS]\cap (SX]{\Big )}{\Bigg 
]}$, we have$$b\le a\le t \mbox { for some } t\in X\cup {\Big( 
}(XS]\cap (SX]{\Big )}.$$If $t\in X$, then $b\le t\in X$, then $b\in 
(X]\subseteq {\Bigg(}X\cup {\Big( }(XS]\cap (SX]{\Big 
)}{\Bigg]}=Q$.\\If $t\in (XS]\cap (SX]$, then $b\le t\in (XS]\cap 
(SX]$, so $$b\in {\Big (}(XS]\cap (SX]{\Big ]}\subseteq {\Bigg(}X\cup 
{\Big (}(XS]\cap (SX]{\Big )}{\Bigg]}=Q.$$[Instead of 2), we could 
also write $((Q]]=(Q]$ (as this holds for any subset of $S$)].

3) If $T$ is a quasi-ideal of $S$ such that $T\supseteq X$, then
$$Q:={\Bigg (}X\cup {\Big( }(XS]\cap (SX]{\Big )}{\Bigg ]}\subseteq 
{\Bigg (}T\cup {\Big( }(TS]\cap (ST]{\Big )}{\Bigg ]}=(T]=T.$$
$\hfill\Box$

The sufficient condition of Proposition 2 in [8] is the following: \\ 
Suppose $X\cap Q\cap Y\subseteq (YQX]$ for every right ideal $X$, 
every left ideal $Y$ and every quasi-ideal $Q$ of $S$. Then $S$ is 
intra-regular.\\Here is its corrected proof: Let $X\subseteq S$. 
Denote by $r(X)$, $l(X)$, $q(X)$ the right ideal, left ideal and the 
quasi-ideal of $S$, respectively, generated by $X$. By hypothesis, we 
have\begin{eqnarray*}X&\subseteq&r(X)\cap q(X)\cap l(X)\subseteq 
{\Big(}l(X)q(X)r(X){\Big]}\\&=&{\Bigg(}(X\cup SX]{\Bigg (}X\cup 
{\Big( }(XS]\cap (SX]{\Big )}{\Bigg ]}(X\cup XS]{\Bigg 
]}\\&=&{\Bigg(}(X\cup SX){\Bigg (}X\cup {\Big( }(XS]\cap (SX]{\Big 
)}{\Bigg )}(X\cup XS){\Bigg ]}\\&\subseteq  &{\Bigg(}(X\cup SX){\bigg 
(}X\cup (XS]{\bigg )} (X\cup XS){\Bigg ]}\\&=&{\Bigg(}{\bigg(}X^2\cup 
SX^2\cup X(XS]\cup SX(XS]{\bigg )}(X\cup XS){\Bigg ]}\end{eqnarray*}
Since $X(XS]\subseteq (X](XS]\subseteq (X^2S]$ and $SX(XS]\subseteq 
(SX](XS]\subseteq (SX^2S]$, we get\begin{eqnarray*}X&\subseteq&
{\Bigg(}{\bigg(}X^2\cup SX^2\cup (X^2S]\cup (SX^2S]{\bigg )}(X\cup 
XS){\Bigg ]}\\&=&{\Big(}X^3\cup SX^3\cup (X^2S]X\cup (SX^2S]X\cup 
X^3S\cup SX^3S\cup (X^2S]XS\cup (SX^2S]XS{\Bigg ]}.\end{eqnarray*}
Since

$(X^2S]X\subseteq (X^2S](X]\subseteq (X^2SX]\subseteq (X^2S]$,

$(SX^2S]X\subseteq (SX^2S](X]\subseteq (SX^2SX]\subseteq (SX^2S]$,

$(X^2S]XS\subseteq (X^2S](XS]\subseteq (X^2SXS]\subseteq (X^2S]$ and

$(SX^2S]XS\subseteq (SX^2S](XS]\subseteq (SX^2SXS]\subseteq (SX^2S]$, 
we obtain$$X\subseteq {\Big(}X^3\cup SX^2S\cup (X^2S]\cup 
(SX^2S]{\Big]}={\Big(}X^3\cup (SX^2S]\cup 
(X^2S]{\Big]}.$$Then\begin{eqnarray*}X^3&\subseteq&{\Big(}X^3\cup 
(SX^2S]\cup (X^2S]{\Big]}X^2\subseteq {\Big(}X^3\cup (SX^2S]\cup 
(X^2S]{\Big]}(X^2]\\&\subseteq&{\Bigg(}{\Big(}X^3\cup (SX^2S]\cup 
(X^2S]{\Big)}X^2{\Bigg]}\\&=&{\Big(}X^5\cup (SX^2S]X^2\cup 
(X^2S]X^2{\Big]}.\end{eqnarray*}Since

$X^5\subseteq SX^2S$

$(SX^2S]X^2\subseteq (SX^2S](X^2]\subseteq (SX^2SX^2]\subseteq 
(SX^2S]$ and

$(X^2S]X^2\subseteq (X^2S](X^2]\subseteq (X^2SX^2]\subseteq (X^2S]$, 
we have$$X^3\subseteq {\Big(}SX^2S\cup (SX^2S]\cup (X^2S]{\Big]}=
{\Big(}(SX^2S]\cup (X^2S]{\Big]}.$$Then\begin{eqnarray*}
X&\subseteq&{\Big(}X^3\cup (SX^2S]\cup (X^2S]{\Big]}\\&\subseteq&
{\Bigg(}{\Big(}(SX^2S]\cup (X^2S]{\Big]}\cup (SX^2S]\cup 
(X^2S]{\Bigg]}\\&=&{\Bigg(}{\Big(}(SX^2S]\cup 
(X^2S]{\Big]}{\Bigg]}\\&=&{\Big(}(SX^2S]\cup 
(X^2S]{\Big]},\end{eqnarray*}and hence\begin{eqnarray*}X^2&\subseteq& 
X{\Big(}(SX^2S]\cup (X^2S]{\Big]}\subseteq (X]{\Big(}(SX^2S]\cup 
(X^2S]{\Big]}\\&\subseteq&{\Bigg(}X{\Big(}(SX^2S]\cup 
(X^2S]{\Big)}{\Bigg]}={\Big(}X(SX^2S]\cup 
X(X^2S]{\Big]}\\&\subseteq&{\Big(}(X](SX^2S]\cup (X](X^2S]{\Big]} 
\subseteq{\Big(}(XSX^2S]\cup 
(X^3S]{\Big]}\\&\subseteq&{\Big(}(SX^2S]{\Big]}=(SX^2S],\end{eqnarray*}

$X^2S\subseteq (SX^2S]S=(SX^2S](S]\subseteq (SX^2S^2]\subseteq 
(SX^2S]$ and

$(X^2S]\subseteq {\Big(}(SX^2S]{\Big]}=(SX^2S]$.\medskip

\noindent Thus we have$$X\subseteq {\Big(}(SX^2S]\cup 
(X^2S]{\Big]}={\Big(}(SX^2S]{\Big]}=(SX^2S],$$ and $S$ is 
intra-regular.$\hfill\Box$\medskip

Combining this result with the Proposition 2 in [8], we get the 
following theorem:\medskip

\noindent{\bf Theorem 1.} {\it Let $S$ be an ordered semigroup. The 
following are equivalent:\begin{enumerate}
\item[$(1)$] S is intra-regular;

\item[$(2)$] For every right ideal X, every left ideal Y and every 
bi-ideal B of S, we have $X\cap B\cap Y\subseteq (YBX]$;

\item[$(3)$] For every right ideal X, every left ideal Y and every 
quasi-ideal Q of S, we have $X\cap Q\cap Y\subseteq 
(YQX]$.\end{enumerate} }\medskip

As we already know, the theory of ordered semigroups based on ideals 
and the theory of $le$-semigroups based on ideal elements are 
parallel to each other. All the results on $le$-semigroups based on 
ideal elements are expressed in ordered semigroups in terms of 
ideals, and conversely. It is surprising that, for the results on 
ordered semigroups based on ideals points do not play any essential 
role but the sets [8], as we have also seen in the results above. In 
this respect, the Theorem 1, in case of
$le$-semigroups is the following:\medskip

\noindent{\bf Theorem 2.} (cf. also [5])  {\it Let S be an 
$le$-semigroup. The following are equivalent:\begin{enumerate}
\item[$(1)$] S is intra-regular;

\item[$(2)$] For every right ideal element x, every left ideal 
element y and every bi-ideal element b of S, we have $x\wedge 
b\wedge y\le ybx$;

\item[$(3)$] For every right ideal element x, every left ideal 
element y and every quasi-ideal element q of S, we have $x\wedge 
q\wedge y\le yqx$.\end{enumerate} }\medskip

Let us first give the necessary definitions, and then we will prove 
the theorem. A $poe$-semigroup is an ordered semigroup (: 
$po$-semigroup) $S$ having a greatest element usually denoted by 
$``e"$ (that is, $e\ge a$ for all $a\in S$) [4]. An $le$-semigroup is 
a $poe$-semigroup which is at the same time a lattice (under the 
operations $\vee$ and $\wedge$) such that $a(b\vee c)=ab\vee ac$ and 
$(a\vee b)c=ac\vee bc$ for all $a,b,c\in S$ [1, 2]. An element $a$ of 
an ordered semigroup $S$ is called a right (resp. left) ideal element 
if $ax\le a$ (resp. $xa\le a$) for all $x\in S$ [1]. If $S$ is a 
$poe$-semigroup, then $a$ is a right (resp. left) ideal element of 
$S$ if and only if $ae\le a$ (resp. $ea\le a$) [4]. An element $a$ of 
a $poe$-semigroup $S$ is called a bi-ideal element of $S$ if $aea\le 
a$, and it is called a quasi-ideal element of $S$ if the element 
$ae\wedge ea$ exists (in $S$) and $ae\wedge ea\le a$ [3]. For an 
$le$-semigroup, we denote by $r(a)$, $l(a)$, $q(a)$ the right ideal 
element, the left ideal element and the quasi-ideal element of $S$, 
respectively, generated by $a$. We have $r(a)=a\vee ae$, $l(a)=a\vee 
ea$. A $poe$-semigroup $S$ is called intra-regular if $a\le ea^2e$ 
for all $a\in S$ (cf., for example [4]). The element $a\vee (ae\wedge 
ea)$ is the quasi-ideal element of $S$ generated by $a$ ($a\in S$). 
In fact:
$${\Big(}a\vee (ae\wedge ea){\Big)}e\wedge e{\Big(a\vee (ae\wedge 
ea)}{\Big)}={\Big(}ae\vee (ae\wedge ea)e{\Big)}\wedge {\Big(}ea\vee 
e(ae\wedge ea){\Big)}.$$Since $ae\wedge ea\le ae$, we have $(ae\wedge 
ea)e\le ae^2\le ae$. Since $ae\wedge ea\le ea$, we have $e(ae\wedge 
ea)\le e^2a\le ea$. Thus we have$${\Big(}a\vee (ae\wedge 
ea){\Big)}e\wedge e{\Big(a\vee (ae\wedge ea)}{\Big)}=ae\wedge ea\le 
a\vee (ae\wedge ea).$$Clearly, $a\vee (ae\wedge ea)\ge a$. If $t$ is 
a quasi-ideal element of $S$ such that $t\ge a$, then $a\vee 
(ae\wedge ea)\le t\vee (te\wedge et)\le t$. \bigskip

\noindent{\bf Proof of the theorem.} $(1)\Longrightarrow (2)$. If $S$ 
is intra-regular, then for any element $a$ of $S$, we have$$a\le 
ea^2e=eaae\le e(ea^2e)(ea^2e)e\le ea^2ea^2e.$$Let now $x$ be a right 
ideal element, $y$ a left ideal element and $b$ a bi-ideal element of 
$S$. Then, for the element $x\wedge b\wedge y$ of $S$, we 
have\begin{eqnarray*}x\wedge b\wedge y&\le& e(x\wedge b\wedge 
y)(x\wedge b\wedge y)e(x\wedge b\wedge y)(x\wedge b\wedge 
y)e\\&\le&(ey)(beb)(xe)\le ybx.
\end{eqnarray*}$(2)\Longrightarrow (3)$. This is because the 
quasi-ideal elements of $S$ are bi-ideal elements of $S$ as well. 
Indeed, if $q$ is a quasi-ideal element of $S$, then $qeq\le qe\wedge 
eq\le q$.\\$(3)\Longrightarrow (1)$. Let $a\in S$. Then $a\le ea^2e$. 
In fact: By hypothesis, we have\begin{eqnarray*}a&\le&r(a)\wedge 
q(a)\wedge l(a)\le l(a)q(a)r(a)\\&=&(a\vee ea){\Big(}a\vee (ae\wedge 
ea){\Big)}(a\vee ae)\\&\le&(a\vee ea)(a\vee ea)(a\vee 
ae)\\&=&(a^2\vee ea^2\vee aea\vee eaea)(a\vee ae)\\&=&a^3\vee 
ea^3\vee aea^2\vee eaea^2\vee a^3e\vee ea^3 e\vee aea^2e\vee 
eaea^2e\\&\le&a^3\vee ea^2e\vee ea^2\end{eqnarray*}Then $$a^3\le 
(a^3\vee ea^2e\vee ea^2)a^2=a^5\vee ea^2ea^2\vee ea^4\le ea^2e.$$Thus 
we have $$a\le ea^2e\vee ea^2, \;a^2\le ea^2ea\vee ea^3\le ea^2e,\; 
ea^2\le ea^2e.$$ Thus we get $a\le ea^2e$, and $S$ is 
intra-regular.$\hfill\Box$ \medskip

\noindent{\bf Remark.} The implication $(1)\Rightarrow (2)$ of 
Theorem 2 holds in a $poe$-semigroup $S$ in general for which for any 
right ideal element $x$, every left ideal element $y$ and every 
bi-ideal element $b$ of $S$, the element $x\wedge b\wedge y$ exists.
\bigskip

\noindent This paper has been submitted to the Bulletin of the 
Malaysian Mathematical Sciences Society on May 2, 2014.

\end{document}